\documentclass[preprint,12pt,authoryear]{elsarticle}

\usepackage{amssymb}

\usepackage{amsfonts,amsmath}
\usepackage{url,hyperref}
\usepackage{booktabs,multirow}
\usepackage{eurosym}
\usepackage{subcaption,epstopdf}
\usepackage{color}

\DeclareMathOperator{\tr}{tr}
\graphicspath{{Figures/}}

\begin{document}

\begin{frontmatter}



\title{Commitment and Dispatch of Heat and Power Units via Affinely Adjustable Robust Optimization}


\author[dtu_compute]{Marco Zugno\corref{cor}}\ead{mazu@dtu.dk}
\author[dtu_compute]{Juan Miguel Morales}\ead{jmmgo@dtu.dk}
\author[dtu_compute]{Henrik Madsen}\ead{hmad@dtu.dk}

\address[dtu_compute]{Deparment of Applied Mathematics and Computer Science, Technical University of Denmark, Matematiktorvet, Building 303B, DK-2800 Kgs. Lyngby, Denmark}
\cortext[cor]{Corresponding author. Tel. +45 4525 3369}

\begin{abstract}
The joint management of heat and power systems is believed to be key to the integration of renewables into energy systems with a large penetration of district heating. Determining the day-ahead unit commitment and production schedules for these systems is an optimization problem subject to uncertainty stemming from the unpredictability of demand and prices for heat and electricity. Furthermore, owing to the dynamic features of production and heat storage units as well as to the length and granularity of the optimization horizon (e.g., one whole day with hourly resolution), this problem is in essence a multi-stage one. We propose a formulation based on robust optimization where recourse decisions are approximated as linear or piecewise-linear functions of the uncertain parameters. This approach allows for a rigorous modeling of the uncertainty in multi-stage decision-making without compromising computational tractability. We perform an extensive numerical study based on data from the Copenhagen area in Denmark, which highlights important features of the proposed model. Firstly, we illustrate commitment and dispatch choices that increase conservativeness in the robust optimization approach. Secondly, we appraise the gain obtained by switching from linear to piecewise-linear decision rules within robust optimization. Furthermore, we give directions for selecting the parameters defining the uncertainty set (size, budget) and assess the resulting trade-off between average profit and conservativeness of the solution. Finally, we perform a thorough comparison with competing models based on deterministic optimization and stochastic programming.

\end{abstract}


\end{frontmatter}

\section*{Nomenclature} \label{sec:nomenclature}

\subsection*{Sets}

\begin{description}
  \item[$t$] time index;
  \item[$k$] unit index;
  \item[$i$] storage index;
\end{description}

\subsection*{Parameters}

\begin{description}
  \item[$\underline q_k$/$\overline q_k$] minimum/maximum heat output of unit $k$;
  \item[$\underline p_k$/$\overline p_k$] minimum/maximum power output of unit $k$;
  \item[$\underline r_k$/$\overline r_k$] downward/upward ramping limit for fuel input of unit $k$;
  \item[$T^\mathrm{U}_k$/$T^\mathrm{D}_k$] minimum up-/down-time for unit $k$;
  \item[$\underline u_i$/$\overline u_i$] minimum extraction/maximum injection of heat from/to storage $i$;
  \item[$\underline s_i$/$\overline s_i$] minimum/maximum heat content of storage $i$;
  \item[$c^0_k$] fixed cost for unit $k$ when on;
  \item[$c_k$] per unit fuel cost of unit $k$;
  \item[$c^\mathrm{SU}_k$/$c^\mathrm{SD}_k$] start-up/shut-down cost for unit $k$;
  \item[$d_t$] heat demand at time $t$.
\end{description}

\subsection*{Decision variables}

\begin{description}
  \item[$v_{kt}$] on/off status of unit $k$ at time $t$ (binary);
  \item[$y_{kt}$] start-up status of unit $k$ at time $t$ (binary);
  \item[$z_{kt}$] shut-down status of unit $k$ at time $t$ (binary);
  \item[$q_{kt}$] heat production from unit $k$ at time $t$;
  \item[$p_{kt}$] power production from unit $k$ at time $t$;\pagebreak[4]
  \item[$u_{it}$] heat injection into storage $i$ at time $t$;
  \item[$s_{it}$] heat content in storage $i$ at time $t$;
\end{description}

\section{Introduction} \label{sec:introduction}

Combined Heat-and-Power (CHP) plants are highly efficient units where exhaust gases from the electricity production process are used to produce heat (co-generation). That way, energy that would otherwise be spilled in the atmosphere is further exploited for residential or industrial heating purposes. It is estimated that CHP units based on the Combined Cycle or Open Cycle Gas Turbine (CCGT, OCGT) technology can reach an efficiency of 90\%, while typical values are 60\% for standard power-only plants based on the same technology \citep{iea11}.

Another important feature of CHP plants is their flexibility. Gas-fired CHP units, just like their power-only counterparts, are able to change the level of their output (ramp-up/-down) relatively quickly. Furthermore, they enjoy relatively short start-up and shut-down time. Moreover, CHP units have two additional levels of flexibility compared to power plants. Firstly, they are typically coupled with heat storage facilities, i.e., insulated tanks where hot water can be stored with negligible losses. This allows for a shift in time between production and consumption of heat. Secondly, the ratio between the heat and power output is variable for the so-called \emph{extraction} CHP plants. This gives them extra flexibility, in that the heat and power outputs can adapt to the demand and price of these commodities.

As a result of their attractive properties in terms of efficiency and flexibility,
co-generation is now experiencing a revival in view of the increasing commitment of industrialized countries to curb CO$_2$ emissions. Indeed, CHP units provide an efficient way of producing heat that can be made essentially carbon neutral by switching from fossil fuel to biomass \citep{iea14}. Furthermore, owing to their flexibility, CHP units can support the integration of large-scale electricity production from variable and partly-predictable renewable sources, e.g., wind and solar, by providing back-up power when the output of the latter is insufficient \citep{iea11}. 
To unlock this potential, it is critical that the heat infrastructure and production assets are thought, planned and managed effectively and jointly with the other energy carriers (electricity, gas), the transport system and the industrial assets in the same urban area \citep{mei13}. This work addresses a fundamental part of this complex problem: the joint optimization of short-term operation of heat and power generation units.

As far as electricity is concerned, our assumption is that the owner of the production assets elects the day-ahead market as the preferred floor for trading their electricity output. This is consistent with the current situation in many regional power markets, where the largest part of the exchange of electricity takes place one day in advance \citep{web10}. As the gate closure for day-ahead markets is typically set around noon on the day before delivery of electricity, participants must schedule electricity production, including the on/off status of the production plants (unit commitment). Similar arrangements requiring the day-ahead scheduling of heat production may also be in place. Indeed, the compatibility of the heat production plan with respect to the transmission and distribution constraints must be verified in advance by the district heating network operator.

As a result of the arrangements described above, the owners of CHP plants must determine the commitment of their units and the production schedule for heat and power with a significant advance in time. Many important parameters of the problem, however, are not known at the day-ahead stage. For example, this is the case of heat demand and electricity prices, both of which are only partly predictable. Another feature of the problem at hand is its multi-stage nature, which results from the fact that the unit commitment and production schedules to be determined span a whole day (typically with a granularity of one hour). This is particularly important as the studied system is inherently dynamic, owing to the presence of intertemporal constraints such as ramping limits and to the state-space representation of storages.

In view of the stochastic and dynamic nature of the problem, the framework of Affinely Adjustable Robust Optimization (AARO) is appropriate, see \cite{ben09}. In our view, this is one of the few mathematical programming techniques that reconcile an appropriate representation of the uncertainty with a rigourous modeling of the multi-stage nature of the system that does not violate the non-anticipativity of the sequence of decisions.

The existing literature on the short-term operation of heat and power systems focused on the use of different optimization techniques. The problem of optimizing the operation of a storage-CHP unit system is tackled via stochastic programming by \cite{pal94} and via stochastic dynamic programming by \cite{rav94}. A deterministic approach for scheduling a district heating system is studied by \cite{ari03}, resulting in a Mixed Integer Linear Programming formulation of the problem. Only more recently did electricity markets enter the picture. \cite{rol04} used a stochastic programming approach to the problem including a model of the uncertainty in heat demand and day-ahead power prices. The stochastic programming approach has been enriched by further considering the intra-day stage by \cite{der14} or by allowing for the submission of offering curves in the day-ahead market by \cite{dim14}.

Recently, a number of papers studied various applications of robust optimization to power systems. The unit commitment problem has been studied via (two-stage) Adaptive Robust Optimization (ARO) considering uncertainty in load \citep{ber12}, wind power production \cite{jia12} or contingencies \cite{str11}. The dispatch problem, where unit commitment variables are not considered, for both energy and reserve has been studied with a similar approach by \cite{zug13} with a view to immunizing the solution to uncertainty in wind power production. \cite{war13} incorporated linear decision rules within a stochastic optimization approach to a similar dispatch problem. Furthermore, \cite{lim15} proposed an ARO model for the short-term optimization of a virtual power plant under uncertainty in prices and wind power production.

In view of the state-of-the-art described above, the contributions of this paper are the following. Firstly, we propose an Affinely Adjustable Robust Optimization (AARO) model for the unit commitment and dispatch problem for units participating in both the electricity and heat markets. Secondly, we set up a case-study using realistic data for unit parameters, heat consumption and prices obtained for the Copenhagen area, which accounts for roughly 20\% of the total consumption of heat in Denmark \citep{var14}. Thirdly, we conduct an extensive numerical study where we assess the use of AARO for the problem at hand from the following perspectives: performance of purely linear decision rules against piecewise linear decision rules, sensitivity of results to changes in the parameters defining the uncertainty set, and comparison with the competing approaches of deterministic optimization and stochastic programming.

The paper is structured as follows. Section \ref{sec:modeling} introduces the deterministic version of the optimization model for unit commitment and dispatch of the heat and power system. Uncertainty in electricity prices and heat demand is then incorporated by reformulating the model as an AARO problem in Section \ref{sec:AARO}. The numerical study based on data from the Copenhagen area is presented and its results are discussed in Section \ref{sec:example}. Finally, conclusions are drawn in Section \ref{sec:conclusion}.

\section{Deterministic modeling of heat and power system} \label{sec:modeling}

This section introduces the deterministic optimization model for the unit commitment and dispatch of a heat and power system. The first part of the section is dedicated to modeling of production units: back-pressure CHP units are dealt with in Section \ref{sec:modeling_backpressure}, extraction CHP plants in Section \ref{sec:modeling_extraction} and heat-only plants in Section \ref{sec:modeling_heat_only}. Then, heat storages are described in Section \ref{sec:modeling_storage}. Finally, Section \ref{sec:modeling_market} joins the unit models in a common market framework. 

\subsection{Back-pressure units} \label{sec:modeling_backpressure}

Back-pressure units are a particular type of CHP units that can only produce heat and power at a fixed ratio. Let $\mathcal K^\mathrm{B}$ be the set of units $k$ of the back-pressure type, and $t$ a time period within the optimization horizon. A back-pressure unit is characterized by a heat-to-power ratio $r^b_k$ such that:
\begin{align}
  p_{kt} = r^b_k \cdot q_{kt} \;, && \forall k \in \mathcal K^{\mathrm B},t \;. \label{eq:backpressure_eq}
\end{align}

The heat output for the unit is limited both from below and from above by $\underline q_k$ and $\overline q_k$, respectively, via the following constraints:
\begin{align}
  \underline q_k v_{kt} \leq q_{kt} \leq \overline q_k v_{kt} \;, && \forall k \in \mathcal K^{\mathrm B},t \;. \label{eq:heat_bounds}
\end{align}
The binary variable $v_{kt} \in \left\{ 0,1 \right\}$ denotes the on/off status of unit $k$, forcing a null heat output when the unit is off.

Given that the power and heat outputs are linearly dependent, there is no need to enforce upper and lower bounds for the power output, which follow trivially from \eqref{eq:backpressure_eq}:
\begin{align}
  \underline p_k = r^b_k \cdot \underline q_k \;, && \forall k \in \mathcal K^{\mathrm B} \;, \\
  \overline p_k = r^b_k \cdot \overline q_k \;, && \forall k \in \mathcal K^{\mathrm B} \;.
\end{align}

We approximate fuel consumption, $f_{kt}$, as a linear function of the heat and power output:
\begin{align}
  f_{kt}  = \varphi_k^\mathrm{p} p_{kt} + \varphi_k^\mathrm{q} q_{kt} && \forall k \in \mathcal K^{\mathrm B},t \;, \label{eq:fuel_consumption}
\end{align}
where $\varphi_k^\mathrm{p}$ and $\varphi_k^\mathrm{q}$ represent the fuel consumption per unit of power and heat output, respectively.

Technical limitations restrain the change in operating regime (ramping) of a unit between consecutive time periods. The following constraints impose lower and upper ramping limits, $\underline r_k$ and $\overline r_k$, on fuel consumption of unit $k$:
\begin{align}
  - \underline r_k \leq f_{kt}  - f_{k(t-1)} \leq \overline r_k && \forall k \in \mathcal K^{\mathrm B},t \;. \label{eq:fuel_ramp}
\end{align}
As an alternative, similar constraints could be imposed on the individual heat and power outputs rather than on the fuel consumption.

The additional binary variables $y_{kt}, z_{kt} \in \left\{ 0,1 \right\}$ indicate, if equal to 1, that unit $k$ is being started up or shut down at time $t$, respectively. They are related mutually and to the unit status, $v_{kt}$, via the following constraints:
\begin{subequations}
  \begin{align}
    & v_{kt} - v_{k(t-1)} - y_{kt} + z_{kt} = 0 \;, && \forall k \in \mathcal K^{\mathrm B},t \;, \label{eq:binary_relations1} \\
    & y_{kt} + z_{kt} \leq 1 \;, && \forall k \in \mathcal K^{\mathrm B},t \;. \label{eq:binary_relations2}
  \end{align}
  \label{eq:binary_relations}
\end{subequations}
Constraint \eqref{eq:binary_relations2} imposes that a unit cannot be started up or shut down at the same time. Note that such a constraint is not needed if the start-up and the shut-down operations have a strictly positive cost.

Minimum up-time is imposed via the following:
\begin{subequations}
  \begin{align}
    & \sum_{t=1}^{T^{\mathrm U 0}_k} v_{kt} = T^{\mathrm U 0}_k \;, && \forall k \in \mathcal K^{\mathrm B} \;, \label{eq:mut0} \\
    & \sum_{\tau=t}^{T^\mathrm{Uf}_k(t)} \left( v_{k\tau} - y_{kt} \right) \geq 0 \;, && \forall k \in \mathcal K^{\mathrm B},t \;. \label{eq:mut_gen}
  \end{align}
  \label{eq:mut}
\end{subequations}
The initial up-time $T^{\mathrm U 0}_k$ is enforced by \eqref{eq:mut0}. Similarly, \eqref{eq:mut_gen} enforces minimum up-time during the other time periods. Note that in the latter, we define
\begin{equation}
  T^\mathrm{Uf}_k(t) = \min{\left\{ t + T^\mathrm{U}_k - 1,T \right\}} \;, \qquad \forall k \in \mathcal K^{\mathrm B},t \;, \label{eq:ut_definition}
\end{equation}
where $T$ indicates the final time period in the horizon and $T^\mathrm{U}_k$ the minimum up-time for unit $k$. The constraints enforcing minimum down-time are specular to \eqref{eq:mut} and skipped here for the sake of conciseness. We refer the reader to \cite{mor14} for further details.


Finally, we define the cost of operating the back-pressure unit $k$ at time $t$ as the sum of fuel, no-load, start-up and shut-down costs:
\begin{align}
  {\rm cost}_{kt}  = c_k f_{kt} + c^0_k v_{kt} + c^\mathrm{SU}_k y_{kt} + c^\mathrm{SD}_k z_{kt} \;, && \forall k \in \mathcal K^{\mathrm B},t \;.
  \label{eq:cost}
\end{align}

\subsection{Extraction units} \label{sec:modeling_extraction}

Extraction units allow for more flexibility in the ratio between heat and power outputs than back-pressure units. Let us consider a unit $k$ of the set $k \in \mathcal K^{\mathrm X}$ of extraction units in the system.

The definition of the feasible region in the heat-power space is based on the following relaxation of the back-pressure constraint \eqref{eq:backpressure_eq}:
\begin{align}
  p_{kt} \geq r^b_k \cdot q_{kt} && \forall k \in \mathcal K^{\mathrm X},t \;. \label{eq:backpressure_ineq}
\end{align}
Like in the case of a back-pressure unit, lower and upper bounds for heat production are defined by \eqref{eq:heat_bounds}.

Fuel consumption is given by \eqref{eq:fuel_consumption}. For an extraction unit, we define upper and lower bounds on fuel consumption:
\begin{align}
  \underline f_k v_{kt} \leq f_{kt} \leq \overline f_k v_{kt} \;, && \forall k \in \mathcal K^{\mathrm X},t \;. \label{eq:fuel_bounds}
\end{align}
Note that, as a result of definition \eqref{eq:fuel_consumption}, \eqref{eq:fuel_bounds} implies that the feasible set in the heat-power space for each unit $k \in \mathcal K^{\mathrm X}$ and at each time period $t$ is included in the region between two parallel lines with the following slope:
\begin{align}
  r^v_k = - \frac{\varphi^q_k}{\varphi^p_k} \;, && \forall k \in \mathcal K^{\mathrm X} \;.
\end{align}
When fuel consumption is at the upper bound, this slope represents the decrease in electricity output needed to produce an extra unit of heat without exceeding the upper bound in \eqref{eq:fuel_bounds}. Hence it is often referred to as \emph{marginal electricity loss for heat production}, see \cite{web05}.

Constraints \eqref{eq:fuel_ramp}--\eqref{eq:mut}, i.e., ramping limits, binary variable constraints, minimum up- and down-time as defined for the back-pressure unit in Section \ref{sec:modeling_backpressure}, hold for an extraction unit too. Similarly, definition \eqref{eq:cost} of operating cost also holds.

\subsection{Heat-only production units} \label{sec:modeling_heat_only}

Let us define unit $k$ of the set $\mathcal K^\mathrm{H}$ as a unit only capable of producing heat. We impose lower and upper bounds on heat production \eqref{eq:heat_bounds}, relations between binary variables \eqref{eq:binary_relations}, minimum up-time \eqref{eq:mut}--\eqref{eq:ut_definition} as well as minimum down-time.

Fuel consumption for these units is solely determined by heat production:
\begin{align}
  f_{kt}  = \varphi^q_k q_{kt} \;, && \forall k \in \mathcal K^\mathrm{H},t \;, \label{eq:fuel_consumption_heat-only}
\end{align}
Fuel-to-heat coefficient $\varphi^q_k$ can in this case be interpreted as the inverse of the unit efficiency.

Ramping limits on fuel consumption \eqref{eq:fuel_ramp} still hold. In this case, they can be directly translated to ramping limits on heat production. Furthermore, the definition of the operating cost \eqref{eq:cost} is also valid for this type of unit.

\subsection{Heat storage units} \label{sec:modeling_storage}

Let us consider a heat storage $i$ from the set $\mathcal I$ of accumulator tanks installed in the system.

The level of heat stored in the tank, $s_{it}$, is governed by the following state-update equation:
\begin{align}
  s_{it} = s_{i(t-1)} + u_{it} \;, && \forall i \in \mathcal{I},t \;,
  \label{eq:state_update}
\end{align}
where $u_{it}$ represents the flow of heat from the production units into the tank at time $t$.

We impose the following lower and upper bounds on the level of heat in the storage:
\begin{align}
  \underline s_i \leq s_{it} \leq \overline s_i \;, && \forall i \in \mathcal{I},t \;.
  \label{eq:state_bounds}
\end{align}
Typically, the lower bound is set to 0, while the upper bound is the thermal capacity of the storage.

Moreover, we assume that the flow of heat to/from the storage tank is bounded from below and above:
\begin{align}
  \underline u_i \leq u_{it} \leq \overline u_i \;, && \forall i \in \mathcal{I},t \;.
  \label{eq:injection_bounds}
\end{align}

Finally, we enforce the condition that the final heat content in the storage be equal to the initial one through:
\begin{align}
  s_{iT} = s_{i0} \;, && \forall i \;.
  \label{eq:state_terminal}
\end{align}
Note that this condition is needed to avoid that the storage be completely emptied throughout the optimization horizon. An alternative solution would consist in assigning an economic value to the heat stored in the tank at the last time period. However, we prefer the solution with terminal condition \eqref{eq:state_terminal} for two reasons. Firstly, the value of stored heat is not easily quantifiable. Secondly, \eqref{eq:state_terminal} makes the comparison between competing optimization methods more transparent, as it ensures that the same amount of heat is stored at the end of the horizon regardless of the modeling choices. In the alternative solution where a value is assigned to the heat stored at the end of the horizon, that comparison would be dependent on the arbitrary choice of the economic value for the terminal heat level.

\subsection{Heat and power market framework} \label{sec:modeling_market}

We consider an electricity pool with time-varying prices set by the interception between supply and demand, see \cite{mor14}. As far as heat is concerned, we consider a market with a fixed price set by a regulator.

We make the assumption that the owner of the CHP system is responsible for satisfying the demand of heat in a given urban area. This is enforced mathematically via the heat balance equation:
\begin{equation}
  \sum_{k \in \mathcal K} q_{kt} - \sum_{i \in \mathcal I} u_{it} = d_t \;, \qquad \forall t\;, \label{eq:heat_balance}
\end{equation}
where $d_t$ is the heat load at time $t$ and $\mathcal K = \mathcal K^\mathrm{B} \cup \mathcal K^\mathrm{X} \cup \mathcal K^\mathrm{H}$.

The objective function is given by the market revenues from electricity sales minus the total costs for operating the units:
\begin{equation}
  \rho = \sum_{t,k \in \mathcal K} \left\{ \lambda_t p_{kt} - c_k \left[ \varphi^p_k p_{kt} +    \varphi^q_k q_{kt} \right] - c^0_k v_{kt} - c^\mathrm{SU}_k y_{kt} - c^\mathrm{SD}_k z_{kt} \right\}  \;,
  \label{eq:obj_function}
\end{equation}
where $\lambda_t$ is the electricity price at time $t$. Note that sales of heat are not included in the objective function, because the price for heat is fixed and so is the production over the optimization horizon, as a result of the combination of \eqref{eq:state_terminal} and \eqref{eq:heat_balance}. Consequently, revenues from heat sales are a constant, and hence they can be omitted from objective function \eqref{eq:obj_function}. 

\section{Decision-making via Affinely Adjustable Robust Optimization} \label{sec:AARO}

The optimization model presented in the previous section is deterministic, in that it relies on the assumption that all parameters be known with certainty at the time of decision-making. However, this assumption is unrealistic for the problem at hand, as unit-commitment and dispatch decisions must be made before the gate closure of the day-ahead markets described in Section \ref{sec:introduction}, i.e., when heat demand and power prices are unknown. In this section, we cast the problem as a two-stage optimization problem within the framework of Affinely Adjustable Robust Optimization (AARO).

\subsection{Formulation as an Affinely Adjustable Robust Optimization problem}

In order to formulate the problem within the AARO framework, we shall first endow the decision-maker with the capability to perform recourse decisions, i.e., adjustable to the realization of the uncertain parameters. Let us collect all the first-stage variables in a vector $\mathbf x$, and define the vector $\mathbf y_{\boldsymbol \delta}$ of recourse decisions as follows:
\begin{align}
  \mathbf x = \begin{bmatrix}
    \mathbf v \\ \mathbf y \\ \mathbf z \\ \mathbf p \\ \mathbf q \\ \mathbf s \\ \mathbf u
  \end{bmatrix} \;, &&
  \mathbf y_{\boldsymbol \delta} = \begin{bmatrix}
    {\boldsymbol \Delta \mathbf p_{\boldsymbol \delta}} \\ {\boldsymbol \Delta \mathbf q_{\boldsymbol \delta}} \\ {\boldsymbol \Delta \mathbf s_{\boldsymbol \delta}} \\ {\boldsymbol \Delta \mathbf u_{\boldsymbol \delta}}
  \end{bmatrix} \;. \label{eq:1st_2nd_var_def}
\end{align}
Unit commitment variables belong to the set of first-stage decisions, as unit status cannot be changed on a real-time basis. The vector $\mathbf x$ further includes the day-ahead dispatch vectors $\mathbf p, \mathbf q, \mathbf s, \mathbf u$. Real-time adjustments with respect to the dispatch are included in the vector of recourse variables $\mathbf y_{\boldsymbol \delta}$.

With the definitions above, we can formulate a stochastic version of the optimization problem described in Section \ref{sec:modeling} as follows:
\begin{subequations} \label{mod:classic}
    \begin{align}
      \underset{\mathbf x,\mathbf y_{\boldsymbol \delta}}{\mathrm{Min.\;}} & \mathbb E_{\boldsymbol \delta} \left\{ \mathbf{c}_{\boldsymbol \delta}^\top \mathbf x + \mathbf{g}_{\boldsymbol \delta}^\top \mathbf y_{\boldsymbol \delta} \right\} \label{eq:classic_obj} \\
      \mathrm{s.t.\;} & \mathbf{A x} \geq \mathbf b \;, \label{eq:classic_1st_const} \\
                      & \mathbf{T^\mathrm{e} x} + \mathbf{W^\mathrm{e} y}_{\boldsymbol \delta} = \mathbf h^\mathrm{e}_{\boldsymbol \delta} \;, \label{eq:classic_link_eq} \\
                      & \mathbf{T^\mathrm{i} x} + \mathbf{W^\mathrm{i} y}_{\boldsymbol \delta} \geq \mathbf h^\mathrm{i}_{\boldsymbol \delta} \;. \label{eq:classic_link_ineq}
    \end{align}
\end{subequations}
We grouped all the constraints involving only first-stage variables into \eqref{eq:classic_1st_const}, while the ones including recourse variables are divided into equalities \eqref{eq:classic_link_eq} and inequalities \eqref{eq:classic_link_ineq}. These last two types of constraints are subject to right-hand side uncertainty stemming from the stochastic nature of the heat demand in \eqref{eq:heat_balance}. The cost coefficients in the objective function in \eqref{eq:classic_obj} are also uncertain as they include electricity market prices, see \eqref{eq:obj_function}. We assume the same marginal cost for dispatch $\mathbf x$ and adjustments $\mathbf y_{\boldsymbol \delta}$. However, the power dispatch is sold at the day-ahead market price, while the real-time adjustment is sold or bought at the balancing market price. We assume price consistency, i.e., that the expectation of the balancing market price be equal to the day-ahead price, see \cite{zav15}.

In order to put us in the framework of AARO, we make the following assumptions:
\begin{enumerate}
  \item[A1.] We require \eqref{eq:classic_link_eq}--\eqref{eq:classic_link_ineq} be valid $\forall \boldsymbol \delta \in \mathcal U$, where the uncertainty set \citep{ber11} $\mathcal U$ is a bounded polyhedron described by the set of $\ell$ linear inequalities $\mathbf L \boldsymbol \delta \geq \mathbf l$.
  \item[A2.] The uncertain parameters $\mathbf g_{\boldsymbol \delta}$, $\mathbf h^\mathrm{e}_{\boldsymbol \delta}$, $\mathbf h^\mathrm{i}_{\boldsymbol \delta}$ depend linearly on a random vector $\boldsymbol \delta$, i.e., $\mathbf g_{\boldsymbol \delta} = \mathbf{\widehat g} + \mathbf G \boldsymbol \delta$, $\mathbf h^\mathrm{e}_{\boldsymbol \delta} = \mathbf{\widehat h^\mathrm{e}} + \mathbf H^\mathrm{e} \boldsymbol \delta$, $\mathbf h^\mathrm{i}_{\boldsymbol \delta} = \mathbf{\widehat h^\mathrm{i}} + \mathbf H^\mathrm{i} \boldsymbol \delta$. Note that there is no loss of generality here, as we could simply redefine the uncertainty vector $\boldsymbol \delta$ as the concatenation of all the stochastic parameters in the model, and also the associated probability space accordingly.
  \item[A3.] The recourse decision $\mathbf y_{\boldsymbol \delta}$ is restricted to be an affine function of the uncertainty, i.e., $\mathbf y_{\boldsymbol \delta} = \mathbf Y \boldsymbol \delta$.
\end{enumerate}

After replacing the affine dependencies to the uncertainty in \eqref{mod:classic}, we obtain the following:
\begin{subequations} \label{mod:RO1}
    \begin{align}
      \underset{\mathbf x,\mathbf Y}{\mathrm{Min.\;}} & \mathbb E_\delta \left\{ \mathbf{c}_{\boldsymbol \delta}^\top \mathbf x + \mathbf{\widehat g}^\top \mathbf Y \boldsymbol \delta + \boldsymbol \delta^\top \mathbf G^\top \mathbf Y \boldsymbol \delta \right\} \label{eq:RO1_obj} \\
      \mathrm{s.t.\;} & \mathbf{A x} \geq \mathbf b \;, \label{eq:RO1_1st_const} \\
                      & \mathbf{T^\mathrm{e} x} + \mathbf{W^\mathrm{e} Y} \boldsymbol \delta = \mathbf{\widehat h^\mathrm{e}} + \mathbf H^\mathrm{e} \boldsymbol \delta \;, & \forall \boldsymbol \delta \in \mathcal U \;, \label{eq:RO1_link_eq} \\
                      & \min_{\boldsymbol \delta \in \mathcal U} \left\{ \left( \mathbf{W^\mathrm{i} Y} - \mathbf H^\mathrm{i} \right) \boldsymbol \delta \right\} \geq \mathbf{\widehat  h^\mathrm{i}} - \mathbf{T^\mathrm{i} x} \;, \label{eq:RO1_link_ineq}
    \end{align}
\end{subequations}
where \eqref{eq:RO1_link_ineq} is equivalent to requiring that \eqref{eq:classic_link_ineq} be valid $\forall \boldsymbol \delta \in \mathcal U$ (note that the $\min$ operator acts row-wise). With some further reformulations we can cast the problem as the following Mixed-Integer Linear Problem (MILP):
\begin{subequations} \label{mod:RO2}
    \begin{align}
      \underset{\mathbf x,\mathbf Y,\boldsymbol \Lambda}{\mathrm{Min.\;}} & \mathbb E_{\boldsymbol \delta} \left\{ \mathbf{c}_{\boldsymbol \delta}^\top \right\} \mathbf x + \mathbf{\widehat  g}^\top \mathbf Y \mathbb E_{\boldsymbol \delta} \left\{ \boldsymbol \delta \right\} + \tr \left\{ \mathbf G^\top \mathbf Y \left( \boldsymbol \Sigma_{\boldsymbol \delta} + \mathbb E \{ \boldsymbol \delta \} \mathbb E \{ \boldsymbol \delta \}^\top \right) \right\} \label{eq:RO2_obj} \\
      \mathrm{s.t.\;} & \mathbf{A x} \geq \mathbf b \;, \label{eq:RO2_1st_const} \\
                      & \mathbf{T^\mathrm{e} x} = \mathbf{\widehat  h^\mathrm{e}} \;, \label{eq:RO2_link_eq1} \\
                      & \mathbf{W^\mathrm{e} Y} = \mathbf H^\mathrm{e} \;, \label{eq:RO2_link_eq2} \\
                      & \boldsymbol \Lambda^\top \mathbf l \geq \mathbf{\widehat h^\mathrm{i}} - \mathbf {T^\mathrm{i} x} \;, \label{eq:RO2_link_const1} \\
                      & \mathbf L^\top \boldsymbol \Lambda = \left( \mathbf{W^\mathrm{i} Y} - \mathbf H^\mathrm{i} \right)^\top \;, \label{eq:RO2_link_const2} \\
                      & \boldsymbol \Lambda \in \mathbb R^{\ell \times m}_{\geq 0} \;, \label{eq:RO2_link_const3} \\
                      & \mathbf Y \in \mathbb R^{n \times r}, Y_{ij} = 0 \;, \forall (i,j) \in \mathcal Z \;. \label{eq:RO2_link_const4}
    \end{align}
\end{subequations}
In order to obtain \eqref{eq:RO2_obj}, we exploited the invariance of the trace operator to cyclic permutations of the operands. The pairs \eqref{eq:RO2_link_eq1}--\eqref{eq:RO2_link_eq2} are equivalent to \eqref{eq:RO1_link_eq}, see \cite{war13}. Note that \eqref{eq:RO2_link_eq1} defines the day-ahead dispatch as the response to the nominal value of the right-hand side $\mathbf{\widehat h^\mathrm{e}}$. Constraints \eqref{eq:RO2_link_const1}--\eqref{eq:RO2_link_const3} are equivalent to \eqref{eq:RO1_link_ineq}. They are obtained by applying strong linear duality row-wise to each of the constraints in \eqref{eq:RO1_link_ineq}, as commonly done in robust optimization to determine the robust counterpart of a single inequality constraint, see \cite{ber11}. The (non-negative) dual variables defined in each of these operations are collected into the $\ell \times m$ matrix $\boldsymbol \Lambda$, where $\ell$ is the number of constraints defining the uncertainty set $\mathcal U$ and $m$ the number of inequality constraints in \eqref{eq:RO1_link_ineq}. Constraints \eqref{eq:RO2_link_const4} define some elements of the matrix $\mathbf Y$ be equal to zero. As explained in Section \ref{sec:AARO_DR}, this is necessary to guarantee the non-anticipativity of the solution, as well as to implement either purely linear or piecewise-linear decision rules.

Note that the objective function in \eqref{eq:RO2_obj} includes the conditional mean of the day-ahead power price (included in $\mathbf c_{\boldsymbol \delta}$), as well as the mean, $\mathbb E \{ \boldsymbol \delta \}$, and the variance-covariance matrix, $\boldsymbol \Sigma_{\boldsymbol \delta}$, of the uncertainty $\boldsymbol \delta$.

\subsection{Definition of budget uncertainty set} \label{sec:AARO_unc_set}

In this section, we briefly introduce the so-called \emph{budget uncertainty set}, which is used throughout this paper to model $\mathcal U$ in the AARO model of the previous section. We refer the reader to \cite{ben09} for a thorough introduction to the concept of uncertainty set.

In this work, we consider budget uncertainty sets defined by the following linear constraints:
\begin{subequations}
  \begin{align}
    & -\mathbf 1 \leq \boldsymbol \delta \leq \mathbf 1 \;, \label{eq:unc_set_interval}\\
    & \left\| \boldsymbol \delta \right\|_1 = \mathbf 1^\top \left| \boldsymbol \delta \right| \leq \Gamma \;. \label{eq:unc_st_budget}
  \end{align}
  \label{eq:unc_set}
\end{subequations}
Inequalities \eqref{eq:unc_set_interval} limit the uncertainty in a symmetric interval centered about zero and of radius equal to one. Note that the affine dependence between the uncertain parameters in optimization model \eqref{mod:classic} and the variables $\boldsymbol \delta$ allows us to consider any nominal values ($\mathbf{\widehat g}$, $\mathbf{\widehat h^\mathrm{e}}$, $\mathbf{\widehat h^\mathrm{i}}$) and interval sizes (which are ultimately set by the matrices $\mathbf G$, $\mathbf H^\mathrm{e}$, $\mathbf H^\mathrm{i}$) for electricity prices and heat demand. Constraint \eqref{eq:unc_st_budget} defines the budget of uncertainty for $\boldsymbol \delta$, enforcing that its total absolute deviation from zero (i.e., summed across all the elements of the vector) is no larger than $\Gamma$. The latter is a parameter that can be tuned by the decision-maker. Larger values of $\Gamma$ imply a more robust solution, as it must be able to cope with uncertainty taking values within a larger set. Note that through the affine dependence of the parameters from $\boldsymbol \delta$, \eqref{eq:unc_st_budget} allows us to enforce rules such as: ``during the considered optimization horizon, an uncertain parameter can have a deviation from the nominal value equal to the interval size during at most $\Gamma$ hourly periods''.

Because of the presence of the absolute value operator, constraint \eqref{eq:unc_st_budget} is nonlinear. However, it can be represented linearly by extending the uncertainty space with the inclusion of positive and negative parts of $\boldsymbol \delta$:
\begin{equation}
\boldsymbol \delta' =
  \begin{bmatrix}
    \boldsymbol \delta \\
    \boldsymbol \delta^+ \\
    \boldsymbol \delta^-
  \end{bmatrix} =
    \begin{bmatrix}
    \boldsymbol \delta \\
    \max(\boldsymbol \delta,\mathbf 0) \\
    - \min(\boldsymbol \delta,\mathbf 0)
  \end{bmatrix} \;,
  \label{eq:unc_lifted}
\end{equation}
where $\max$ and $\min$ operate row-wise. A linear representation of the uncertainty set can be then obtained by replacing \eqref{eq:unc_set} with the following linear constraints:
\begin{subequations}
  \begin{align}
    & \boldsymbol \delta = \boldsymbol \delta^+ - \boldsymbol \delta^- \;, \\
    & \mathbf 0 \leq \boldsymbol \delta^+, \boldsymbol \delta^- \leq \mathbf 1 \;, \\
    & \mathbf 1^\top \left( \boldsymbol \delta^+ + \boldsymbol \delta^- \right) \leq \Gamma \;.
  \end{align}
  \label{eq:unc_set_linear}
\end{subequations}

\subsection{Definition of decision rules} \label{sec:AARO_DR}

Decision rules are defined by the matrix $\mathbf Y$, which is a decision variable in optimization model \eqref{mod:RO2}. Matrix $\mathbf Y$ defines the real-time adjustment of units and storage operation as a function of the realization of the uncertainty $\boldsymbol \delta$. In principle, the definition of $\mathbf Y$ in \eqref{mod:RO2} may establish links between adjustments of operation variables and uncertainties whose value has not been revealed yet, hence destroying the \emph{non-anticipativity} of the decision structure. Thus, additional constraints must be imposed to enforce that the elements of $\mathbf Y$ linking the readjustment to uncertainty unfolding at future time periods be zero. If there is one element of $\boldsymbol \delta$ per time period and the vector is ordered chronologically, this corresponds to requiring that $\mathbf Y$ is the concatenation of lower triangular matrices (one matrix per type of readjustment variable).

Moreover, let us consider the extended definition \eqref{eq:unc_lifted} of the uncertainty space including positive and negative variables. Similarly, we can define extended decision rules $\mathbf Y'$ including, besides the $\mathbf Y$ coefficients for $\boldsymbol \delta$, also $\mathbf Y^+$ and $\mathbf Y^-$ for $\boldsymbol \delta^+$ and $\boldsymbol \delta^-$, respectively. We can define purely linear decision rules by enforcing $\mathbf Y^+, \mathbf Y^- = \mathbf 0$. Similarly, piecewise-linear decision rules are defined by setting $\mathbf Y = \mathbf 0$.

Constraint \eqref{eq:RO2_link_const4} enforces that some elements of the $\mathbf Y$ matrix be zero under a choice of the index set $\mathcal Z$ that is consistent with the two observations above, i.e., it complies with the non-anticipativity of the solution and properly implements either linear or piecewise-linear decision rules.

We conclude this section with an important remark.
When using piecewise linear decision rules, after replacing $\mathbf Y$ with $\mathbf Y'$ and $\boldsymbol \delta$ with $\boldsymbol \delta'$ in the objective function in \eqref{eq:RO2_obj}, we obtain (omitting zero terms):
\begin{multline}
  \mathbb E_{\boldsymbol \delta} \left\{ \mathbf{c}_{\boldsymbol \delta}^\top \right\} \mathbf x + \mathbf{\widehat  g}^\top \mathbf Y^+ \mathbb E_{\boldsymbol \delta} \left\{ \boldsymbol \delta^+ \right\} + \mathbf{\widehat  g}^\top \mathbf Y^- \mathbb E_{\boldsymbol \delta} \left\{ \boldsymbol \delta^- \right\} \\
  + \tr \left\{ \mathbf G^\top \mathbf Y^+ \mathbb E_{\boldsymbol \delta} \{ \boldsymbol \delta^+ \boldsymbol \delta^\top \} + \mathbf G^\top \mathbf Y^- \mathbb E_{\boldsymbol \delta} \{ \boldsymbol \delta^- \boldsymbol \delta^\top \} \right\} \;. \label{eq:RO2_obj_PLDR}
\end{multline}
Hence, the use of piecewise linear decision rules calls for a richer model of the uncertainty than with purely linear decision rules, as it requires the estimation of means and expectations of products between $\boldsymbol \delta$, $\boldsymbol \delta^+$ and $\boldsymbol \delta^-$.

\section{Case study} \label{sec:example}

In this section, we present and discuss results obtained from an extensive numerical study assessing various features of the proposed model.

We consider a system comprising an extraction and a back-pressure CHP unit, a peaker producing only heat and a heat storage. The technical parameters for the production units are listed in Table \ref{tab:unit_parameters}. Parameters for the CHP units are based on values for existing units of the Copenhagen area as reported in \cite{vph14}. Notably, the back-pressure unit is the cheapest. However, it is less flexible than the extraction unit, which, besides having larger ramping limits and shorter minimum up- and down-time, can produce heat and power at different ratios. The heat-only unit is meant to be employed only as a back-up, i.e., during periods of peak demand. Hence, we assumed larger costs for operating this unit. This unit is also very flexible in terms of minimum and maximum heat output, ramping and minimum up- and down-times. We make the assumption that the back-pressure unit cannot modify the production schedule at the real-time stage, hence it is marked as ``non-flexible'' in Table \ref{tab:unit_parameters}. This models a unit that for technical reasons (e.g., the type of fuel used) is too slow to perform changes in its output level with short notice. We consider a storage tank with a heat capacity $\overline s = 2000\,$MWh and maximum in-/out-flow $\overline u = - \underline u = 300\,$MW.
\begin{table}[htb]
  \centering
  \caption{Parameters for the units in the illustrative example}
  \label{tab:unit_parameters}
  \begin{tabular}{lccrrr}
      \toprule
      \multicolumn{2}{l}{\multirow{2}{*}{Description and symbol}} & \multirow{2}{*}{Unit} & \multicolumn{2}{c}{CHP} & \multirow{2}{*}{Heat-only} \\
      \cmidrule{4-5}
                                            &             &       & extr.         & back-pr.      & \\
      \midrule
      Power/heat ratio                      & $r_b$             & --        & 0.64      & 0.28      & 0 \\
      Power loss for heat prod.             & $r_v$             & --        & $-$0.13     & --        & -- \\
      Fuel per electricity unit             & $\varphi^p$       & --        & 2.40      & 2.40      & -- \\
      Fuel per heat unit                    & $\varphi^q$       & --        & 0.31      & 0.36      & 1.09 \\
      Min. fuel input                       & $\underline f$    & MWh       & 120       & 72.24     & 0 \\
      Max. fuel input                       & $\overline f$     & MWh       & 631.20    & 516       & 1086.96 \\
      Ramp-up limit                         & $\overline r$     & MWh/h     & 150       & 50        & 1086.96 \\
      Ramp-down limit                       & $\underline r$    & MWh/h     & 150       & 50        & 1086.96 \\
      Min. heat output                      & $\underline q$    & MWh       & 0         & 70        & 0 \\
      Max. heat output                      & $\overline q$     & MWh       & 331       & 500       & 1000 \\
      Min. power output                     & $\underline p$    & MWh       & 41.56     & 19.60     & 0 \\
      Max. power output                     & $\overline p$     & MWh       & 263       & 140       & 0 \\
      Fuel cost                             & $c$               & \euro/MWh & 24.16     & 12.75     & 93.96 \\
      No-load cost                          & $c^0$             & \euro     & 0         & 0         & 2684.56 \\
      Start-up cost                         & $c^\mathrm{SU}$   & \euro     & 7\,382.55 & 6\,040.27 & 0 \\
      Shut-down cost                        & $c^\mathrm{SD}$   & \euro     & 7\,382.55 & 6\,040.27 & 0 \\
      Min. up-time                          & $T^\mathrm{U}$    & h         & 2         & 5         & 0 \\
      Min. down-time                        & $T^\mathrm{D}$    & h         & 2         & 5         & 0 \\
      Flexible                              & --                & --        & yes       & no        & yes \\
      \bottomrule
  \end{tabular}
\end{table}

We employ real-world data for electricity prices and heat demand in the Copenhagen region in Denmark. Power prices for Eastern Denmark (DK2 area of Nord Pool) are publicly available for download at \cite{ene15}. We use heat consumption data from the western Copenhagen area, VEKS, available at \cite{mad15}. Note that the latter data refer to the time period between July 1995 and June 1996. Since the power market was not liberalized at that point in time, we paired heat load data with power price data from the corresponding week in the period 2013--2014. We rescaled the heat consumption data linearly so that the heat-only plant covers roughly 0.5\% of the total load over the year. We consider four weeks representative of summer (the 4th week after the start of the dataset), autumn (17th), winter (30th) and spring (43rd). The winter week includes the largest hourly realization of heat load over the entire year.

In the remainder of this section, we consider the following methods:
\begin{enumerate}
  \item Deterministic Optimization (DET): corresponding to the model described in Section \ref{sec:modeling} with expected values in place of $\lambda_t$ and $d_t$;
  \item Robust Optimization with Linear Decision Rules (RO-LDR): implementing \eqref{mod:RO2} using purely linear decision rules;
  \item Robust Optimization with Piecewise-Linear Decision Rules (RO-PLDR): implementing \eqref{mod:RO2} using piecewise-linear decision rules;
  \item Stochastic Programming (SP): a discretization of \eqref{mod:classic} with sampling of the uncertainty $\mathbf \delta$ resulting in a finite number of scenarios.
\end{enumerate}
In the implementation of the robust optimization models, we considered that decision rules only depend on the uncertainty in heat demand and not on the realization of the balancing market price. This reduces the modeling of the uncertainty to the estimation of price and heat load expectations as well as their variance-covariance matrix. Without loss of generality, we made the assumption of no autocorrelation of the processes and slightly positive correlation (0.3) between heat load and balancing market prices. The variances we used for the heat demand and balancing market prices are consistent with the state-of-the-art on forecasting, which reports RMSE values of about 7\% for heat load \citep{nie06,zel13} and 33\% for balancing market prices \citep{jon14}. In the case of piecewise linear decision rules, we calculated the mean of the positive and negative parts as:
\begin{equation}
  \mathbb E \left\{ \delta^+ \right\} = \mathbb E \left\{ \delta^- \right\} = \frac{\mathbb E \left\{ \left| \delta \right| \right\}}{2} = \sigma_{\delta} \sqrt{\frac{1}{2 \pi}} \;,
\end{equation}
which holds for a zero-mean Normal variable, see \cite{leo61}. We further assume that positive and negative parts of the heat load forecasting error contribute equally to the correlation with the balancing market price.

As far as the stochastic programming model is concerned, we generated scenarios by drawing samples from a multivariate Normal distribution fully defined by the characteristics described above. We employed 100 scenarios after reduction from an initial set of 2000 performed using the fast-forward method presented in \cite{hei03}. We employ scenario fans, i.e., a collection of individual paths for the uncertainty branching out at the root node (representing the first-stage) but with no further branches at subsequent stages. Note that the non-anticipativity of the solution is not fulfilled within this setup. This approximation is needed as the use of scenario trees within a multi-period problem with 24 stages is prohibitive.

The remainder of this section is structured as follows. In Section \ref{sec:example_operation}, we highlight operational differences between the proposed robust optimization model and its deterministic counterpart. Section \ref{sec:example_LDR_vs_PLDR} assesses the improvement obtained by switching from linear to piecewise linear decision rules within the robust optimization framework. Section \ref{sec:example_price_robustness} presents a sensitivity analysis showing how the parametrization of the uncertainty set affects the performance and the robustness of the solution. Finally in Section \ref{sec:example_RO_vs_SOTA}, we compare the proposed approach with deterministic optimization and stochastic programming.

\subsection{Daily operation} \label{sec:example_operation}

\begin{figure}
  \centering
  \includegraphics{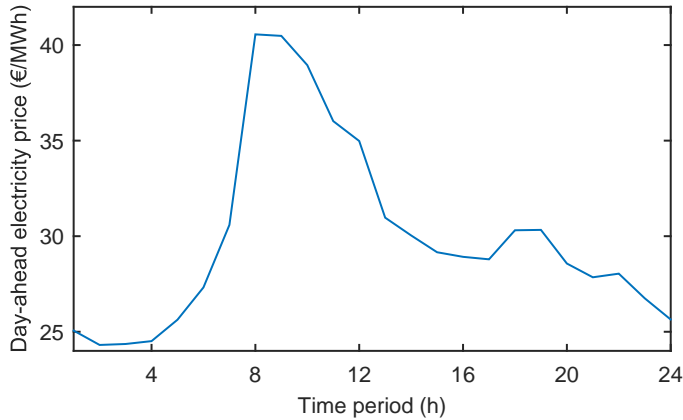}
  \caption{Day-ahead electricity price during simulation day in spring}
  \label{fig:price_w43}
\end{figure}

Let us consider the operation of the heat and power system during the fourth day of the spring week considered in this numerical study. The day-ahead electricity price for this day is shown in Figure \ref{fig:price_w43}. Notably, this price peaks during the morning hours.

\begin{figure}
  \centering
  \begin{subfigure}[t]{\textwidth}
    \includegraphics[width=\textwidth]{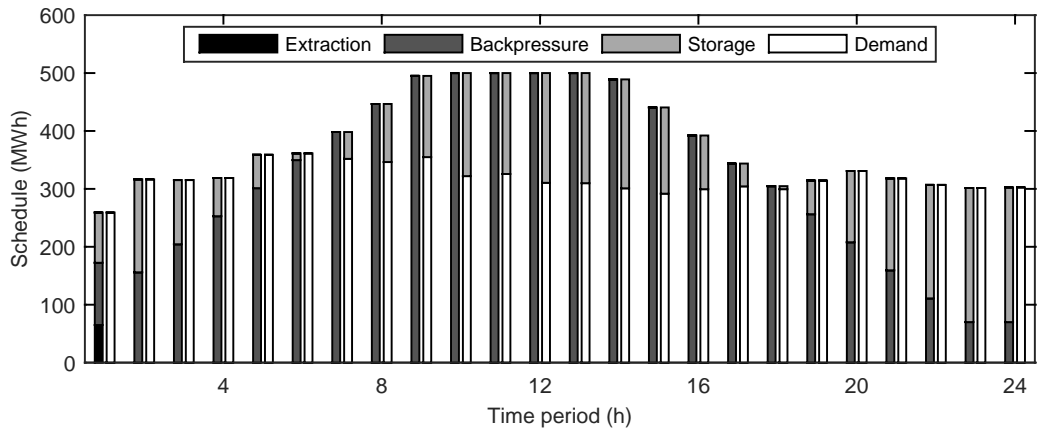}
    \caption{Deterministic model}
    \label{fig:sched_DET_w43}
  \end{subfigure}
  \begin{subfigure}[t]{\textwidth}
    \includegraphics[width=\textwidth]{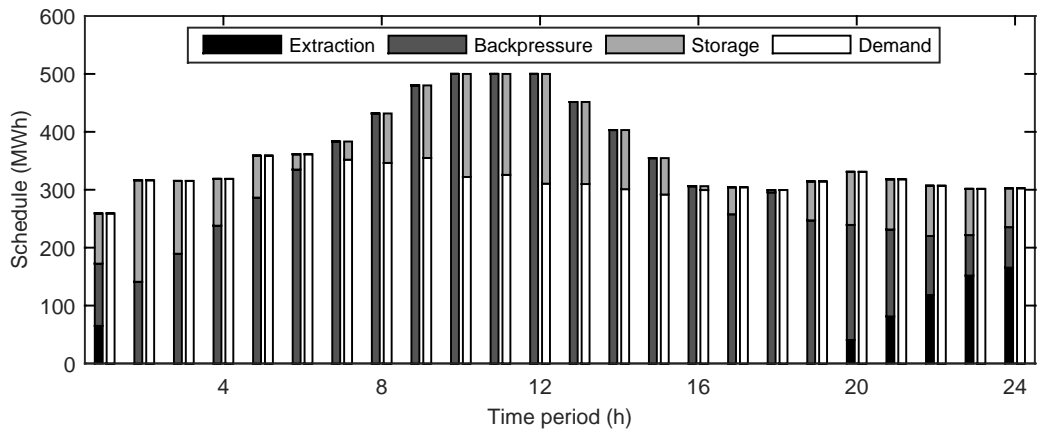}
    \caption{RO-LDR (interval radius equal to 3.2 times the standard deviation, $\Gamma = 6$)}
    \label{fig:sched_RO_w43}
  \end{subfigure}
  \caption{Day-ahead dispatch for the system during simulation day in spring}
  \label{fig:sched_w43}
\end{figure}

The heat production and consumption schedule resulting from the deterministic model is illustrated in Figure \ref{fig:sched_DET_w43}. For each hour, the left bar represents supply and the right one demand of heat. Note that the heat storage can be on both sides depending on whether it is being discharged or charged. Notably, as the heat load is relatively low in this season, the back-pressure CHP unit is sufficient for covering the demand alone (the extraction unit is initially on and quickly turned off). The heat storage is used to perform arbitrage first by being charged as the back-pressure unit ramps up heat and power production in the morning, thus exploiting the electricity price peak, and then by being discharged during the lower power price periods.

The schedule obtained for the same day with the robust optimization model with linear decision rules (RO-LDR) is shown in Figure \ref{fig:sched_RO_w43}, where we employed an interval radius equal to 3.2 times the estimated standard deviation of the heat demand at each time period and a budget of uncertainty $\Gamma = 6$. Two observations can be drawn from the comparison with the deterministic schedule in Figure \ref{fig:sched_DET_w43}. Firstly, the extra conservativeness of the robust optimization approach reduces the extent to which arbitrage is performed: notice, for example, the shorter peak of electricity production around 11am. The second observation concerns the scheduling of the extraction CHP unit at the end of the day. Although more expensive than the back-pressure unit, this plant is needed to be online for the system to cope with possible deviations of the heat demand from the forecast value.

We conclude this section with Table \ref{tab:XT_dispatch_Gamma}, which shows the total planned production from the extraction unit over the 24-hour period as a function of the budget of uncertainty $\Gamma$. Notably, as the requirements in terms of conservativeness of the solution increase, the higher the schedule for the extraction unit. When scheduled, this unit can both ramp up and down in response to unexpected deviations of the heat load.
\begin{table}
  \centering
  \caption{Change of dispatch for extraction CHP unit during simulation day in spring as a function of $\Gamma$}
  \label{tab:XT_dispatch_Gamma}
  \begin{tabular}{lcccccc}
  \toprule
  & \multirow{2}{*}{Unit} & \multicolumn{5}{c}{$\Gamma$} \\
  \cmidrule{3-7}
  & & 2 & 4 & 6 & 8 & 10 \\
  \midrule
  Total heat dispatch & MWh & 265.20 & 444.27 & 622.24 & 811.31 & 1250.74 \\
  Periods online & h & 4 & 5 & 6 & 8 & 24 \\
  \bottomrule
  \end{tabular}
\end{table}

\subsection{Linear vs piecewise linear decision rules} \label{sec:example_LDR_vs_PLDR}

In the previous section, the robust optimization model with linear decision rules is seen to partly replace in the dispatch the back-pressure CHP unit with the extraction one, despite the higher cost of the latter, in order to guarantee flexibility.

Table \ref{tab:XT_dispatch_Gamma_DR} shows that the dispatch decision suggested by the robust optimization model with purely linear decision rules can be suboptimal. This table illustrates the total heat dispatch for the extraction unit over the fourth day of the winter week. Because the peak heat load occurs during this day, the system needs both the back-pressure and the extraction CHP plants to be online, and must use the peaker to provide the remainder of the heat load. Yet, the dispatch of the extraction unit decreases as the uncertainty budget $\Gamma$ grows. This behavior is caused by the linearity of the decision rule, which requires the extraction unit to provide both upward and downward regulation in response to heat load deviations. Hence, the setpoint for heat production from this unit must be lower than the production capacity if this unit has to respond to deviations in heat demand. Piecewise-linear decision rules allow for two different re-dispatch strategies depending on whether the heat load is lower than expected (in which case the extraction CHP unit reduces production) or higher (the peak unit increases production). With purely linear decision rules, these two re-dispatch strategies must be one and the same (the extraction CHP unit is used to cover both positive and negative deviations of the heat load).
\begin{table}
  \centering
  \caption{Change of dispatch for extraction CHP unit during simulation day in winter as a function of $\Gamma$ for linear and piecewise-linear decision rules}
  \label{tab:XT_dispatch_Gamma_DR}
  \begin{tabular}{lcccccc}
  \toprule
  \multirow{2}{*}{Decision rule} & \multirow{2}{*}{Unit} & \multicolumn{5}{c}{$\Gamma$} \\
  \cmidrule{3-7}
  & & 2 & 4 & 6 & 8 & 10 \\
  \midrule
  Linear & MWh & 7944.00 & 7783.77 & 7572.89 & 7377.16 & 7200.06 \\
  Piecewise-linear & MWh & 7944.00 & 7944.00 & 7944.00 & 7944.00 & 7944.00 \\
  \bottomrule
  \end{tabular}
\end{table}

\begin{figure}
  \centering
  \begin{subfigure}{0.48\textwidth}
    \includegraphics[width=\columnwidth]{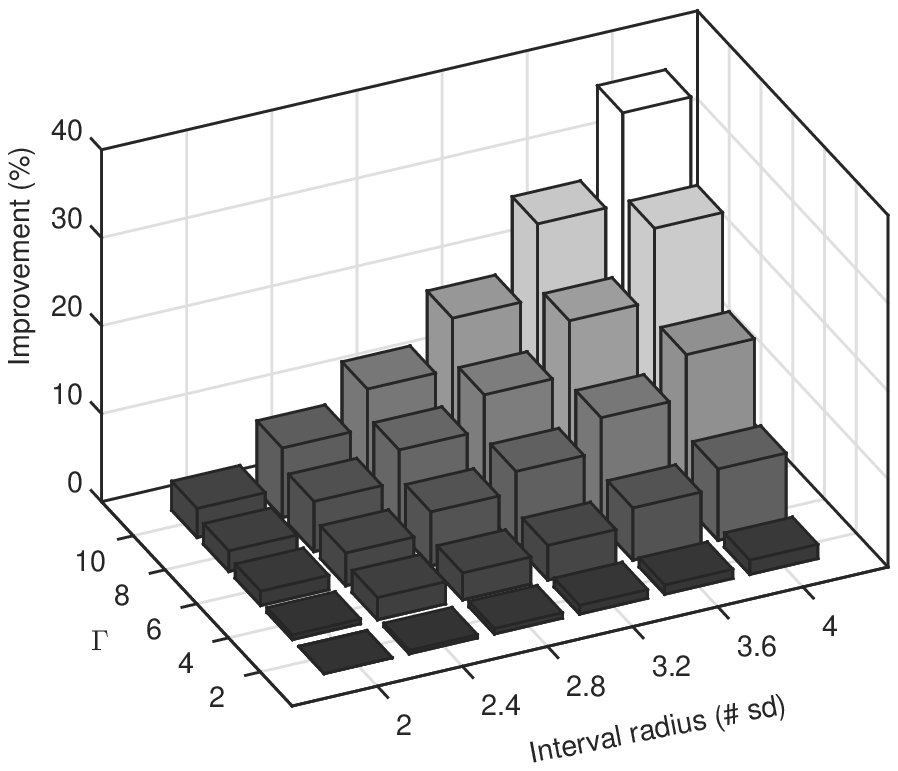}
    \caption{Theoretical}
    \label{fig:LDR_vs_PLDR_in}
  \end{subfigure}
  \hfill
  \begin{subfigure}{0.48\textwidth}
    \includegraphics[width=\columnwidth]{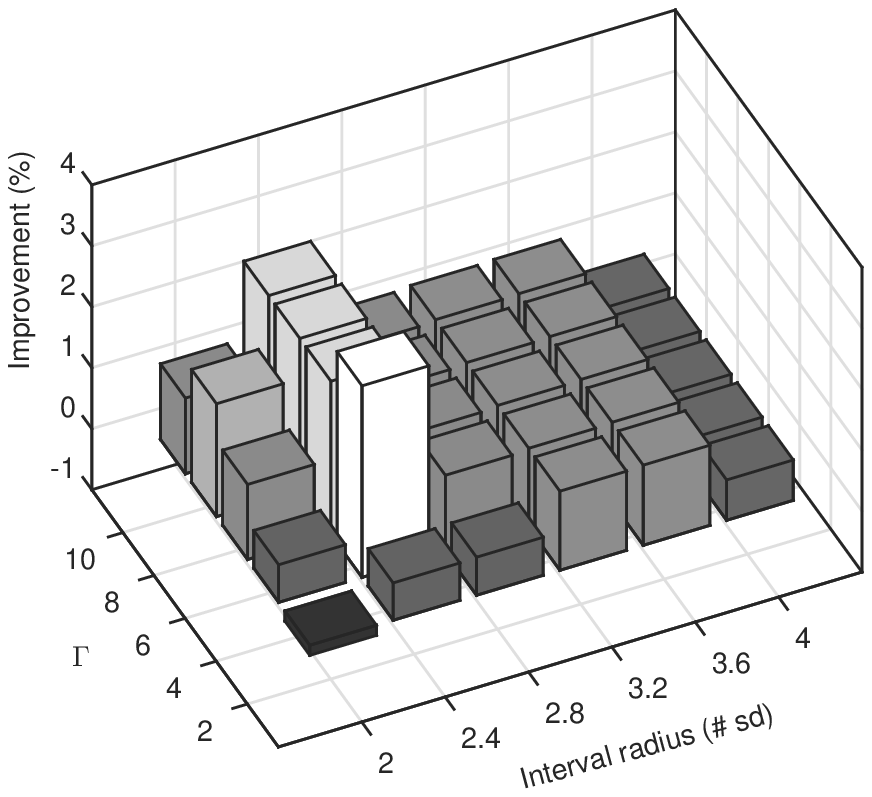}
    \caption{Out-of-sample}
    \label{fig:LDR_vs_PLDR_out}
  \end{subfigure}
  \caption{Improvement in expected profit when switching from linear to piecewise-linear decision rules on winter day}
\end{figure}

The suboptimality of the dispatch described above has a significant impact in terms of profit for the system. Figure \ref{fig:LDR_vs_PLDR_in} shows the improvement in the objective function value obtained when switching from linear to piecewise-linear decision rules. Notably, this theoretical improvement is rather large and increases with the required conservativeness of the solution, both in terms of size of the interval and budget $\Gamma$ for the uncertainty set defined in Section \ref{sec:AARO_unc_set}.
The improvement just mentioned, though, is only theoretical as it is calculated on the assumption that the decision-maker will follow linear or piecewise-linear decision rules in practice, which is not necessarily true. An approximation of the practical improvement can be determined by fixing the unit commitment and dispatch solution of the robust optimization model in a stochastic programming model implementing a discretized version of \eqref{mod:classic}. The solution is then evaluated out-of-sample by considering the redispatch (not necessarily affine on the heat load deviation) over a set of 100 scenarios. Average results for the same day are shown in Figure \ref{fig:LDR_vs_PLDR_out}. Notably, the improvement is lower but still significant. However, it should be remarked that the stochastic programming model does not preserve the non-anticipativity of the solution, as discussed at the beginning of this section. Hence, this also represents an approximation of the actual improvement of piecewise-linear decision rules.

\subsection{Assessing the price of robustness} \label{sec:example_price_robustness}

One aspect of critical importance in robust optimization is fine-tuning the parameters governing the size of the uncertainty set. A certain level of robustness is required by the very nature of the problem at hand. Indeed, failure to meet heat demand may have some important consequences, either financial or in terms of corporate image. On the other hand, increasing robustness requirements are met at the expense of the overall financial performance (in expectation).

\begin{figure}
  \centering
  \begin{subfigure}{0.48\textwidth}
    \includegraphics[width=\columnwidth]{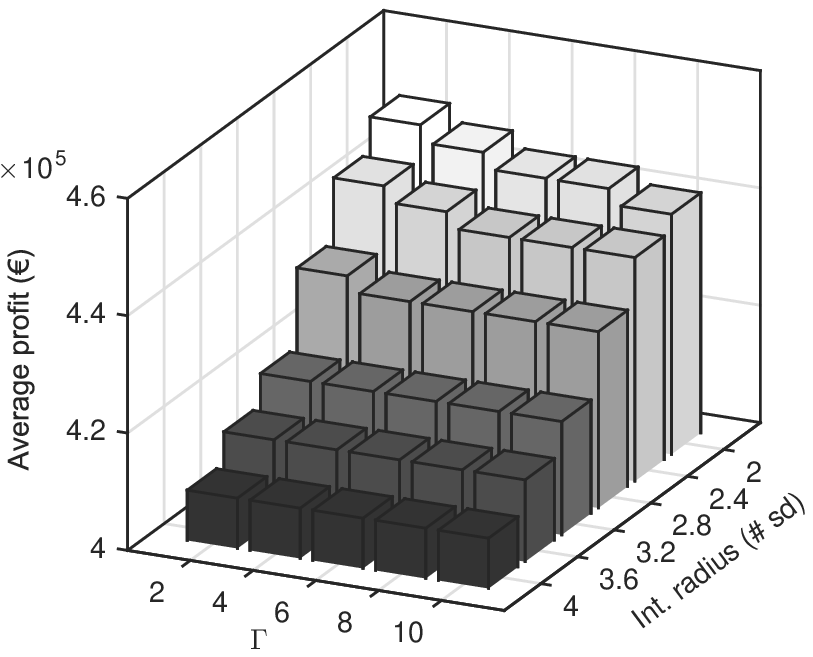}
    \caption{Average profit}
    \label{fig:sensitivity_RO_w30_prof}
  \end{subfigure}
  \hfill
  \begin{subfigure}{0.48\textwidth}
    \includegraphics[width=\columnwidth]{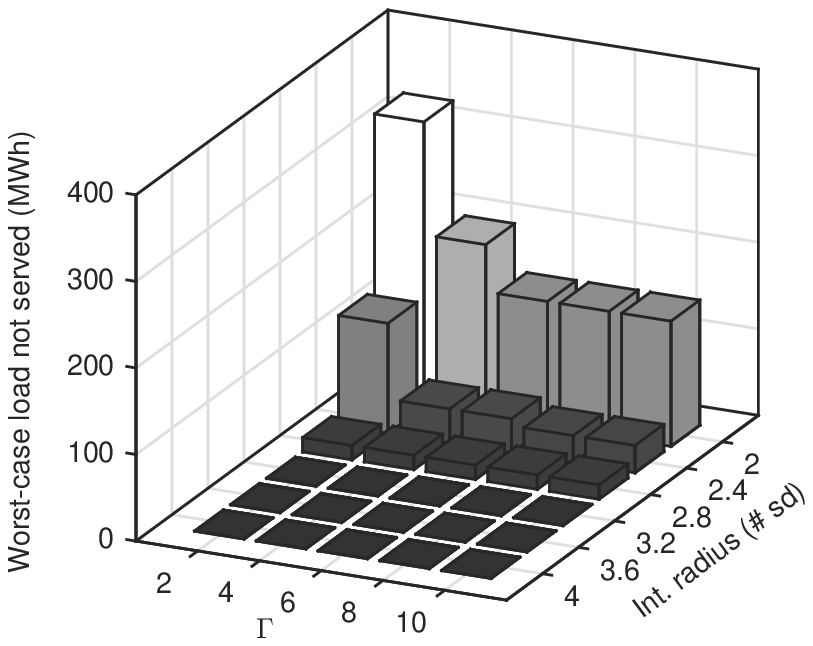}
    \caption{Heat-load not served (worst case)}
    \label{fig:sensitivity_RO_w30_worstLNS}
  \end{subfigure}
  \caption{Out-of-sample average profit and heat-load not served obtained with the robust optimization model under different parameters on winter day}
  \label{fig:sensitivity_RO_w30}
\end{figure}

Figure \ref{fig:sensitivity_RO_w30} illustrates the trade-off between average profit and the worst-case realization of heat-load not served during the winter day characterized by peak heat load. The results in the figure are obtained from out-of-sample simulation where the unit commitment and dispatch solution is fixed in the stochastic programming model. Notably, the average profit decreases as the conservativeness requirement (i.e., interval size and uncertainty budget) is strengthened in Figure \ref{fig:sensitivity_RO_w30_prof}. However, increasing the robustness requirement is beneficial only up to a certain point in terms of heat-load not served. Indeed as Figure \ref{fig:sensitivity_RO_w30_worstLNS} confirms, there appears to be no need to push the interval radius of the uncertainty set much beyond three times the standard deviation.

For the sake of brevity, we omit similar plots obtained for spring, summer and autumn, where we observe similar trends with in general lower requirements in terms of conservativeness. Note that those results are partly included in the following section.

\subsection{Robust optimization vs alternative methods} \label{sec:example_RO_vs_SOTA}

This section is dedicated to a comparison of the robust optimization approach with two other state-of-the-art techniques described at the beginning of this section: deterministic optimization and stochastic programming. We evaluate the unit commitment and dispatch solutions from the different models by fixing them in a stochastic programming approach and evaluating their results after re-sampling the uncertainty (i.e., with different scenarios than the ones used in the optimization of the stochastic programming model).

Results in terms of average profit and heat-load not served (abbreviated as ``heat LNS'') from the winter, spring, summer and autumn week are shown in Tables \ref{tab:paretoWinter}, \ref{tab:paretoSpring}, \ref{tab:paretoSummer} and \ref{tab:paretoAutumn}, respectively. For the robust optimization approach we only included parameter choices that result in Pareto-optimal points in the space of expected profit and worst-case heat-load not served. The average profit does not include penalties for shedding heat load, nor it discounts the fuel saved by the operator of the system as a result of curtailments.

Both the robust optimization and the stochastic programming approach, as expected, significantly reduce the heat-load not served obtained with the deterministic optimization model. Furthermore, if the average heat-load not served is discounted from the average profit with the sum of avoided fuel cost and a penalty (or the fuel cost of an expensive backup unit), both models outperform the deterministic approach. This is particularly true in winter and spring, where the expected heat load not served is rather high for the deterministic solution, see Tables \ref{tab:paretoWinter} and \ref{tab:paretoSpring}.

Moreover, the robust optimization approach with piecewise-linear decision rules provides a flexible framework where conservativeness of unit commitment and dispatch can easily be traded-off with average financial performance simply by tuning the interval size and the budget of the uncertainty set. For sufficiently large values of these parameters, heat-load shedding can be completely removed in all the cases considered. This is particularly important in comparison to the stochastic programming approach, which incurs some load shedding in all weeks excluding the summer one. We remark that this heat-load shedding is caused by the necessary atomization of the sample space of the uncertainty, despite the use of a large sample set with a scenario reduction technique. We deem this drawback rather critical, as it would not be possible to get rid of the residual load shedding by considering a risk measure in the stochastic programming formulation.

The results in Tables \ref{tab:paretoWinter}--\ref{tab:paretoAutumn} show that the stochastic programming solution achieves higher average profits for comparable levels of heat-load not served. This behavior can be explained at least partly by the fact that stochastic programming is used to evaluate all the models. Indeed, this evaluation scheme does not reward the non-anticipativity of the solution, which robust optimization models with decision rules guarantee (on the contrary of stochastic programming), and thus it penalizes the former.

\begin{table}
  \centering
  \caption{Out-of-sample simulation results on winter day. Only Pareto-optimal versions of the robust optimization approach are shown}
  \label{tab:paretoWinter}
  \begin{tabular}{lrrrrr}
    \toprule
    \multicolumn{3}{c}{Model specification} & \multirow{2}{*}{Avg. profit (\euro)} & \multicolumn{2}{c}{Heat LNS (MWh)} \\
    \cmidrule{1-3} \cmidrule{5-6}
    method & radius (\# sd) & $\Gamma$ (--) && largest & expected \\
    \midrule
    \multirow{8}{*}{RO-PLDR} & 2.00 & 2 & 449\,301.37 & 329.38 & 7.88 \\
    & 2.00 & 4 & 446\,337.45 & 198.93 & 3.38 \\
    & 2.00 & 6 & 443\,598.68 & 144.68 & 2.51 \\
    & 2.40 & 2 & 443\,474.44 & 127.20 & 0.51 \\
    & 2.40 & 4 & 440\,776.45 & 39.64 & 0.29 \\
    & 2.40 & 8, 10 & 438\,086.22 & 31.58 & 0.20 \\
    & 2.80 & 2 & 432\,707.89 & 16.59 & 0.05 \\
    \vspace{3pt} & 3.20 & 2--10 & 419\,281.79 & 0.00 & 0.00 \\
    \vspace{3pt} DET & -- & -- & 455\,815.78 & 745.78 & 69.52 \\
    SP & -- & -- & 448\,945.73 & 159.90 & 2.13 \\
    \bottomrule
  \end{tabular}
\end{table}

\begin{table}
  \centering
  \caption{Out-of-sample simulation results on spring day. Only Pareto-optimal versions of the robust optimization approach are shown}
  \label{tab:paretoSpring}
  \begin{tabular}{lrrrrr}
    \toprule
    \multicolumn{3}{c}{Model specification} & \multirow{2}{*}{Avg. profit (\euro)} & \multicolumn{2}{c}{Heat LNS (MWh)} \\
    \cmidrule{1-3} \cmidrule{5-6}
    method & radius (\# sd) & $\Gamma$ (--) && largest & expected \\
    \midrule
    \multirow{3}{*}{RO-PLDR} & 2.00 & 2 & 251\,514.22 & 217.26 & 5.70 \\
    & 2.40 & 2 & 250\,505.55 & 6.78 & 0.01 \\
    \vspace{3pt} & 4.00 & 2 & 248\,484.94 & 0.00 & 0.00 \\
    \vspace{3pt} DET & -- & -- & 260\,927.88 & 340.15 & 49.74 \\
    SP & -- & -- & 257\,295.35 & 76.17 & 1.54 \\
    \bottomrule
  \end{tabular}
\end{table}

\begin{table}
  \centering
  \caption{Out-of-sample simulation results on summer day. Only Pareto-optimal versions of the robust optimization approach are shown}
  \label{tab:paretoSummer}
  \begin{tabular}{lrrrrr}
    \toprule
    \multicolumn{3}{c}{Model specification} & \multirow{2}{*}{Avg. profit (\euro)} & \multicolumn{2}{c}{Heat LNS (MWh)} \\
    \cmidrule{1-3} \cmidrule{5-6}
    method & radius (\# sd) & $\Gamma$ (--) && largest & expected \\
    \midrule
    RO-PLDR & 2.00 & 2 & 88\,423.27 & 0.00 & 0.00 \\
    DET & -- & -- & 94\,413.41 & 117.00 & 16.10 \\
    SP & -- & -- & 92\,356.87 & 0.00 & 0.00 \\
    \bottomrule
  \end{tabular}
\end{table}

\begin{table}
  \centering
  \caption{Out-of-sample simulation results on autumn day. Only Pareto-optimal versions of the robust optimization approach are shown}
  \label{tab:paretoAutumn}
  \begin{tabular}{lrrrrr}
    \toprule
    \multicolumn{3}{c}{Model specification} & \multirow{2}{*}{Avg. profit (\euro)} & \multicolumn{2}{c}{Heat LNS (MWh)} \\
    \cmidrule{1-3} \cmidrule{5-6}
    method & radius (\# sd) & $\Gamma$ (--) && largest & expected \\
    \midrule
    \multirow{2}{*}{RO-PLDR} & 2.00 & 2.00 & 249\,631.81 & 90.35 & 0.41 \\
    \vspace{3pt} & 2.40 & 2.00 & 248\,673.51 & 0.00 & 0.00 \\
    DET & -- & -- & 259\,135.23 & 238.72 & 12.01 \\
    SP & -- & -- & 255\,836.17 & 99.71 & 0.54 \\
    \bottomrule
  \end{tabular}
\end{table} 

\section{Conclusions} \label{sec:conclusion}

In this paper, we address the determination of the optimal day-ahead unit commitment and dispatch of systems including heat and power production units as well as heat storages. Owners of such units face this problem as they have to decide their trading strategy on the day-ahead power market and communicate their schedules for heat production in advance. As a result of the look-ahead time, this optimization problem is subject to the uncertainty in heat demand and in electricity prices, as these are unknown at the time of decision-making. Furthermore, the problem is inherently dynamic (and multi-stage) because of the intertemporal constraints governing the production and storage units.

We propose a formulation of the problem based on robust optimization where optimal recourse decisions are approximated by linear or piecewise-linear functions of the uncertain parameters. Uncertainty is accounted for in this formulation through the definition of an \emph{uncertainty set} as well as of (conditional) means and variance-covariance matrices for the uncertain parameters, i.e., heat demand and power prices. Furthermore, affine decision rules allow the modeler to enforce dependence of each recourse decision solely on the subset of stochastic parameters whose actual value unfolds before the implementation of the recourse decision itself. In other words, the proposed model complies with the non-anticipativity condition in multi-stage decision-making processes.

An extensive numerical study based on data from the Copenhagen area in Denmark highlights several properties of the proposed model. Through an example of daily operation, we show that the higher requirements in terms of conservativeness of the solution (i.e., its capability to cope with unforeseen heat-load deviations) in the proposed model slightly reduce the possibility of exploiting the flexibility of the heat system to perform arbitrage in the electricity market. Further financial results from the study suggest that a significant improvement can be obtained by switching from linear to piecewise-linear decision rules. Moreover, we perform a sensitivity study assessing the trade-off between average profit and conservativeness in terms of heat-load not served resulting from tuning the parameters (size and budget) defining the uncertainty set. Finally, we compare the proposed model with the alternative approaches of deterministic optimization and stochastic programming. We show that, differently from the competing approaches, our model (if appropriately tuned) can result in complete immunization from the heat-load uncertainty. Under a mild penalization of heat-load not served (e.g., including the activation cost of backup units), the average profit can improve that obtained with deterministic optimization. If solely considering average profit, stochastic programming seems to outperform our model. However, this evaluation is performed via stochastic programming (after resampling the uncertain parameters from the assumed distribution), which does not enforce non-anticipativity of the recourse decisions. Devising (and testing the models on) a non-anticipative evaluation scheme is an interesting direction for future research.

Other interesting topics for future research include considering piecewise-linear decision rules with more than two pieces, the inclusion of binary variables in the recourse decisions, the use of more sophisticated models of the uncertainty (e.g., based on time-series models) and the inclusion of more market stages in the model (e.g., reserve or intra-day electricity markets).


\section*{Acknowledgements}

The work of the authors is partly funded by DSF (Det Strategiske Forskningsr{\aa}d) through the \href{http://smart-cities-centre.org/}{CITIES} (no. 1035-00027B) and the \href{http://www.ensymora.dk/}{ENSYMORA} (no. 10-093904) projects, as well as by the \href{http://www.ipower-net.dk/}{iPower} platform project, supported by DSF and RTI (R{\aa}det for Teknologi og Innovation).



\bibliographystyle{elsarticle-harv}
\bibliography{biblio}


\end{document}